\documentclass{article}

\usepackage{amsfonts,amsmath}
\usepackage{amssymb,latexsym}
\usepackage[dvips]{epsfig}
\usepackage{amscd}  
\usepackage{pst-plot}

\title{Patterns in knot cohomology I}
\author{Mikhail Khovanov} 
\date{January 30, 2002}

\newtheorem{prop}{Proposition}
\newtheorem{theorem}{Theorem}
\newtheorem{lemma}{Lemma}
\newtheorem{corollary}{Corollary}
\newtheorem{conjecture}{Conjecture} 

\begin{document}
\maketitle
\baselineskip 12pt

\def\Z{\mathbb Z}
\def\Q{\mathbb Q}
\def\l{\lbrace}
\def\r{\rbrace}
\def\o{\otimes}
\def\lra{\longrightarrow}
\def\mc{\mathcal} 
\def\drawing#1{\begin{center} \epsfig{file=#1} \end{center}}
\def\cH{\mc{H}}
\def\cA{\mc{A}}
\def\cC{\mc{C}}

\begin{abstract}
Cohomology theory of links, introduced in [Kh1], is combinatorial. 
Dror Bar-Natan recently wrote a program that found ranks 
of cohomology groups of all prime knots with up to 11 crossings [BN]. 
His surprising experimental data is discussed in this note. 
\end{abstract}

\tableofcontents

\section{Notations}

\begin{figure} [htb] \drawing{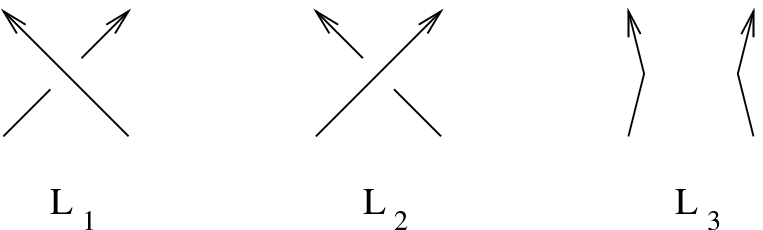} \caption{} \label{skein} 
\end{figure}

The Jones polynomial is determined by the skein relation 
\begin{equation*}
  q^2 J(L_1) - q^{-2} J(L_2) = (q - q^{-1}) J(L_3),
\end{equation*}
where $L_i$ are depicted in figure~\ref{skein}, and by the
normalization  $J(\mathrm{unknot}) = 1.$ 
This standard normalization is different from the one in [Kh1,2].  

Familiarity with [Kh1,2] or [BN] is assumed. \emph{Warning:} we use the 
grading conventions of [Kh2], and  the cohomology group that we denote 
by $\mc{H}^{i,j}$ is denoted $\mc{H}^{i,-j}$ in [BN] and [Kh1].
Let $h^{i,j}(K)$ (or simply $h^{i,j}$) be the rank of $\cH^{i,j}(K).$ 
Ranks of cohomology groups satisfy (notice $q^{-j}$, rather than $q^j$) 
\begin{equation} \label{orig-charac} 
  (q+q^{-1})J(K) = \sum_{i,j} (-1)^i q^{-j} h^{i,j}(K).   
\end{equation}

We use the Rolfsen enumeration for knots with 10 or fewer crossings. 
Knots with more than 10 crossings are enumerated as in
\emph{Knotscape}, for instance, $11^n_{77}$ denotes the $77$th 
non-alternating 11-crossing knot.

\section{Initial observations}\label{observations}

\begin{figure} \drawing{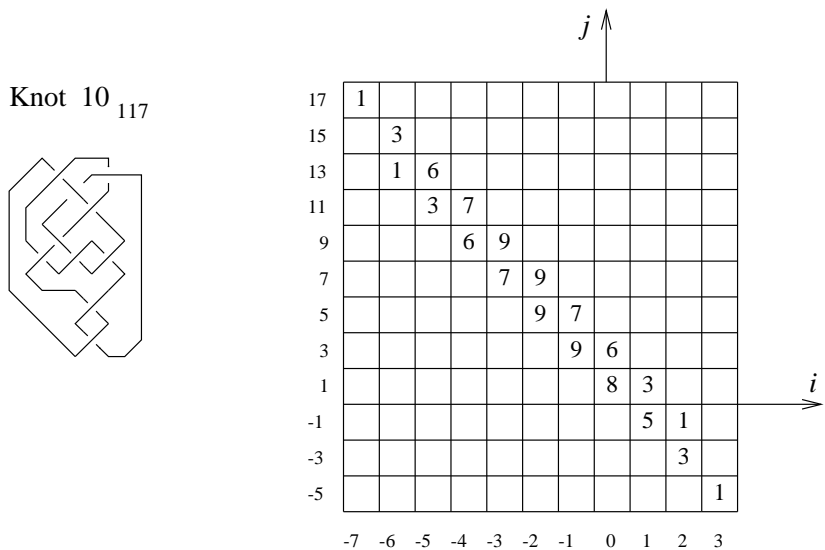} 
\caption{$10_{117}$ and ranks of its cohomology groups} 
\label{knot117} 
\end{figure}

There are 249 prime unoriented knots with  at most 10 crossings 
(not counting mirror images). From Bar-Natan [BN] we learn that 
for all but 12 of these
knots the nontrivial cohomology groups lie on two adjacent diagonals.  
Let us call such knots \emph{homologically thin}, or H-thin, for short. 
We have no clue why nearly all small knots are H-thin. 
Figure~\ref{knot117} depicts $10_{117},$  an H-thin knot, and ranks of its 
cohomology groups. $h^{i,j}$ is zero if the $(i,j)$-square 
is empty. 

Squares with even $j$-coordinates are omitted from the picture, since 
cohomology groups $\cH^{i,2k}(K),$ for a knot $K,$ are always zero. By a 
\emph{diagonal} we mean a line $2i+j=b,$ for some $b,$ also referred
to as the $b$-diagonal. 

All H-thin knots with up to 10 crossings share the following properties 

(i) cohomology groups are supported on $(\sigma\pm 1)$-diagonals, 
where $\sigma$ is the signature of the knot; 

(ii) after substracting $1$ from $h^{0, \sigma \pm 1},$ 
the numbers on the upper diagonal coincide with numbers on 
the lower diagonal after the $(1,-4)$ shift;

(iii) the Jones polynomial is alternating:  
$J(K) = \sum c_i q^{2i},$ if $c_ic_j>0$ then $j \equiv i (\mbox{mod
}2),$ if $c_ic_j<0$ then $j \not\equiv i (\mbox{mod }2).$ 
Unless the knot is a $(2,n)$-torus 
knot, for $n\in \{3,5,7,9\},$ the Jones polynomial has no gaps, i.e. 
$c_i\not= 0, c_{i+k}\not= 0$ implies $c_{i+m}\not= 0$ for all 
$m$ between $1$ and $k-1.$ 

(iv) The Alexander polynomial $\Delta(K)= \sum a_i t^i$ is alternating
and has no gaps.  

\vspace{0.1in}

All alternating and the majority of non-alternating knots with up 
to 10 crossings are H-thin. Knots that are not H-thin will be 
called H-thick (homologically thick). 
The twelve H-thick knots with at most 10 crossings are 

  \begin{center} $8_{19}, 9_{42}, 10_{124}, 10_{128}, 10_{132}, 
     10_{136}, 10_{139}, 10_{145}, 10_{152}, 10_{153}, 10_{154}, 10_{161}.$ 
  \end{center}

\begin{figure} \drawing{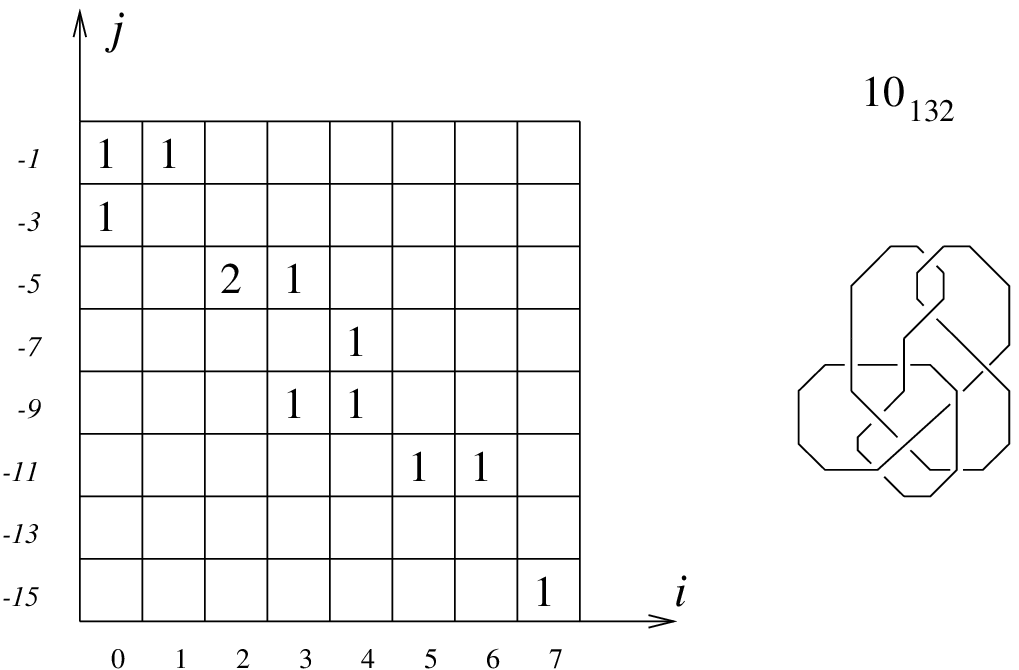} \caption{$10_{132}$ and ranks of
its cohomology groups} \label{knot132} 
\end{figure}

Figure~\ref{knot132} shows the knot $10_{132}$ and ranks of its 
cohomology groups. 

Properties (i), (iii),and (iv) of H-thin knots (with at most 10 crossings) 
fail on many of these knots. The 12 H-thick knots satisfy 

(i') cohomology groups are supported on three adjacent diagonals. 
Discard the diagonal with the smallest total rank of cohomology groups 
supported on it. The two remaining ones are $(\sigma \pm 1)$-diagonals. 

\begin{figure} \drawing{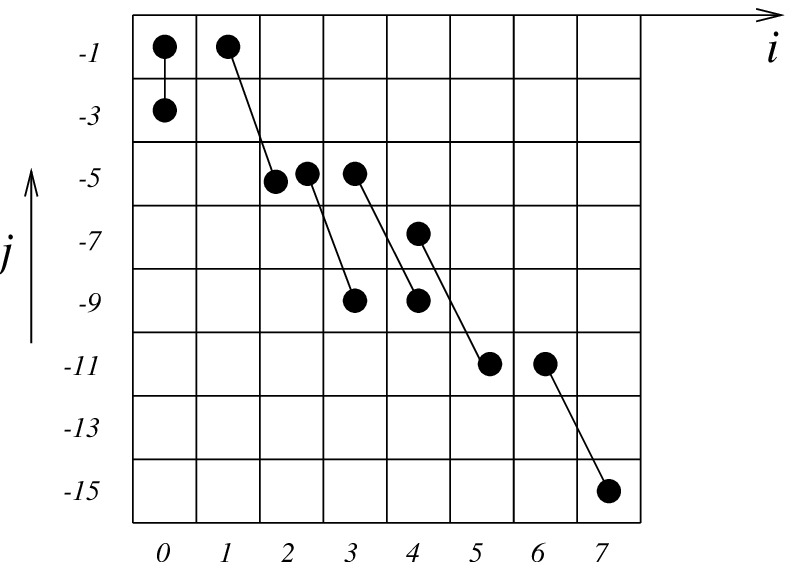} \caption{Cohomology of $10_{132}$
 arranged in pairs} \label{pair132} 
\end{figure}

(ii') if, for a suitable $i,$ we substract $1$ from $h^{0,i}$ and 
$h^{0,i+2},$ the remaining numbers can be arranged into  pairs with 
the $(1,-4)$ difference in the bigrading (figure~\ref{pair132} does it 
for $10_{132}$).

(iii') The Jones polynomials of $10_{124},$  $10_{139},$  $10_{145},$
 $10_{152},$  $10_{153},$  $10_{154},$  $10_{161}$ 
are not alternating. 
The Jones polynomials of $8_{19},$  $10_{124},$  $10_{132},$ $10_{139},$ 
 $10_{145},$ $10_{152},$ $10_{153},$ $10_{154},$ $10_{161}$ have gaps. 

(iv') The Alexander polynomials of $8_{19},$ $10_{124},$ $10_{128},$ 
 $10_{139},$ $10_{145},$ $10_{152},$ $10_{153},$ $10_{154},$ $10_{161}$ 
are not alternating. The Alexander polynomials of $8_{19},$
 $10_{124},$ $10_{139},$ $10_{154},$   $10_{161}$ have gaps. 

We verified (iii), (iii'), (iv), and (iv') using the tables in [St]. 

\vspace{0.2in}

For any knot $K$ the Alexander polynomial at $-1$ equals 
 the Jones polynomial at $\sqrt{-1}:$ 
  \begin{equation} 
  \Delta_{-1}(K) =  J_{\sqrt{-1}}(K)
 \end{equation}
(because of our choice of variable $q,$ the R.H.S. is
 $J_{\sqrt{-1}}(K)$ rather than the more common $J_{-1}(K)$). 

Coefficient-wise, with notations from (iii),(iv), 
  \begin{equation*}
    \sum_i (-1)^i a_i = \sum_i (-1)^i c_i. 
  \end{equation*}
Since Jones and Alexander polynomials of H-thin 
knots with at most 10 crossings are alternating, for these knots we
  obtain  
  \begin{equation*}
     \sum_i | a_i| = \sum_i |c_i|. 
  \end{equation*}
Properties (i) and (ii) imply that, in addition, 
  \begin{equation*}
     \mathrm{rank}\cH(K) -1 = \sum_i |c_i|
  \end{equation*}
where $\mathrm{rank}\cH(K) = \sum_{i,j} h^{i,j}$ is the rank of 
 total cohomology of the knot. To summarize, H-thin knots with at 
most 10 crossings satisfy 
  \begin{equation}
    \mathrm{rank}\cH(K) -1 = \sum_i |c_i| = \sum_i |a_i|. 
  \end{equation}

What about the twelve H-thick knots? For each of them the inequalities 
hold 
  \begin{equation} \label{unequal}
    \sum_i |a_i| \le \mathrm{rank}\cH(K) -3 \ge \sum_i |c_i|,
  \end{equation} 
(Note that $\sum_i |a_i|$ and $\sum_i |c_i|$ are odd for any knot and 
 $\mathrm{rank}(\cH)$ is even.) 

\vspace{0.2in}

Alternating knots with at most 10 crossings are H-thin, 
and it was conjectured in [BN] and [G] that all alternating knots are 
H-thin. This conjecture is now a theorem, due to Eun Soo Lee [Lee]: 

\begin{theorem} Non-split alternating links are H-thin. 
\end{theorem}

\vspace{0.1in} 

We now look at the data for 11-crossing knots [BN]. 
There are 367 alternating and 185 non-alternating prime knots with 
11 crossings. H-thick knots among them number 41. 

Properties (i)-(iv) continue to hold for H-thin knots with 11
crossings. There are several 11-crossing knots with non-alternating 
Jones or Alexander polynomial. All of them are H-thick. Likewise, 
11-crossing knots with a gap in the Alexander polynomial are H-thick.   

\vspace{0.1in} 

{\bf Problem:} \emph{Explain why so many non-alternating knots with 
11 or fewer crossings are H-thin.}

\section{$A$-module structure of knot cohomology}
\label{mod-structure}

In this section we work over $\Q$ 
(rather than over $\Z,$ as in [Kh1, Section 7]). 
In particular, the base ring is $\cA= \Q [X]/(X^2)$ and the chain 
complex $\cC(D)$ associated to a plane diagram $D$ of a knot $K$ is a 
complex of $\Q$-vector spaces. Cohomology groups $\cH^{i,j}(D)$ are
finite-dimensional $\Q$-vector spaces, only finitely many of them 
are nontrivial. Dimensions $h^{i,j}$ of these groups are invariants of $K.$ 

\begin{figure} [htb] \drawing{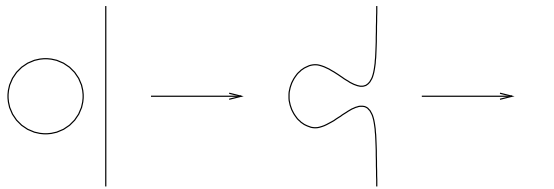} \caption{Cobordism between 
  $\mathrm{circle}\cup D$ and $D$} \label{segment} 
\end{figure}

Choose a segment $I$ of $D$ that does not contain crossings. 
Place an unknotted 
circle next to $I$ and consider the cobordism that merges the circle 
and $I$ (figure~\ref{segment}).
This cobordism induces a map of complexes 
$\cA\otimes \cC(D) \to \cC(D)$ and makes $\cC(D)$ into a complex of graded 
$\cA$-modules. A Reidemeister move from $D$ to $D'$ that happens away 
from $I$ induces a chain homotopy equivalence between complexes 
of $\cA$-modules $\cC(D)$ and $\cC(D').$ Given two diagrams $D_1$ and $D_2$ 
of $K$ and two segments $I_1$ and $I_2$ in them, there is a 
sequence of Reidemeister moves that takes $(D_1,I_1)$ to $(D_2,I_2)$ 
such that all moves happens away from $I_1.$ Instead of moving 
an arc over or under $I_1$ we can move it across the rest of the 
plane (or $\mathbb{S}^2$). In other words, 
there are as many knots as one-component $(1,1)$-tangles.

We obtain an invariant of $K,$ the complex $\cC(D)$ of free $A$-modules up 
to chain homotopy equivalence. The Krull-Schmidt theorem, valid for
bounded complexes of finite-dimensional modules over
finite-dimensional algebras, tells us that $\cC(D)$ decomposes (uniquely
up to an isomorphism) as 
direct sum of an acyclic complex and indecomposable complexes with 
nontrivial cohomology. 
The multiplicity of each indecomposable complex in this decomposition 
is an invariant of $K.$ Denote by $C_n$ the complex  
\begin{equation}\label{indec}
  0 \lra \cA \stackrel{X}{\lra} \cA\{ - 2\} \stackrel{X}{\lra} \cdots 
  \stackrel{X}{\lra} \cA\{ -2n+2 \} \stackrel{X}{\lra} \cA\{ - 2n \} \lra 0, 
\end{equation}
where the leftmost $\cA$ is in cohomological degree $0.$ 

\begin{prop} A non-acyclic indecomposable complex of free graded $A$-modules 
is isomorphic to $C_n[i]\{j\}$ for a unique triple  $(n,i,j), n\ge 0.$ 
\end{prop} 

\emph{Example:} If $K$ is a $(2,2m+1)$-torus knot, $C(K)$ is a 
direct sum of $C_0\{ 2m\}$ and $C_1[2i+1] \{4i+2m+2\}, 1\le i \le m.$

\begin{prop} \label{tensor-prod} 
   $\cC(K_1\# K_2) \cong \cC(K_1) \otimes_{\cA} \cC(K_2).$
\end{prop} 

\emph{Proof:} Obvious. $\square$ 

\vspace{0.1in}

Define \emph{homological width} of $K,$ denoted $hw(K),$ 
as the minimal number $m$ such that 
cohomology of $K$ lie on $m$ adjacent diagonals. The homological 
width of a knot is at least 2, since cohomology groups of  
indecomposable complexes $C_n$ lie on 2 adjacent 
diagonals, and any knot has nontrivial cohomology (since the Jones
polynomial does not vanish). According to our definitions, a knot 
is H-thin if and only if it has homological width 2.  

\vspace{0.1in}

Proposition~\ref{tensor-prod} implies 

\begin{prop} $hw(K_1\# K_2) = hw(K_1) + hw(K_2) - 2. $
\end{prop} 

\begin{corollary} \label{con-sum} $K_1 \# K_2$ is $H$-thin if and only
    if both $K_1$ and $K_2$ are $H$-thin. 
\end{corollary}

\vspace{0.1in}

{\bf Reduced cohomology} 

\vspace{0.1in}

Let $Q= \cA/X\cA$ be the one-dimensional representation of $\cA.$ 
Define the \emph{reduced complex} of $D$ by 
\begin{equation*}
      \widetilde{\cC}(D) = \cC(D) \otimes_{\cA} Q. 
\end{equation*}
This is a complex of graded $\Q$-vector spaces. We call its 
cohomology the \emph{reduced cohomology} of $D$ (and $K$) and 
denote by $\widetilde{\cH}(D)$ and $\widetilde{\cH}(K),$ the latter 
are defined up to isomorphism. Ranks of 
cohomology groups $\widetilde{\cH}^{i,j}(D)$ are invariants of $K.$
The Euler characteristic of $\widetilde{\cH}$ 
is the Jones polynomial (compare to (\ref{orig-charac})): 
\begin{equation*} 
  J(K)  = \sum_{i,j}(-1)^i q^{-j} \hspace{0.02in}
  \mathrm{rank}(\widetilde{\cH}^{i,j}(K)), 
\end{equation*} 
 therefore,  
\begin{equation*}
 \mathrm{rank}\widetilde{\cH}(K) \ge |J_{\sqrt{-1}}(K)| = |
 \Delta_{-1}(K)|. 
\end{equation*}  

\begin{prop} Reduced cohomology groups $\widetilde{\cH}^{i,j}(K)$ lie
on one diagonal ($2i+j$ is constant) if and only if  $K$ is $H$-thin. 
\end{prop} 

\begin{corollary} The Jones polynomial of an $H$-thin knot is
alternating. The absolute values of its coefficients are dimensions 
of reduced cohomology groups. 
\end{corollary}

\vspace{0.1in}

{\bf H-restricted knots}

\vspace{0.1in}

Properties (ii),(ii') admit a homological interpretation. 
We say that a knot $K$ is \emph{$H$-restricted} if non-acyclic
indecomposable summands of the $\mc{A}$-module complex 
$\mc{C}(K)$ are one $\mc{A}\{ i\},$ for 
some $i,$ and one or several $C_1[j]\{ k\},$ for $j,k\in \Z.$ 
Cohomology groups of a $H$-restricted knot can be paired up as in 
(ii'). Existence of such pairing, however, does not imply that 
a knot is $H$-restricted. 

$(2,2m+1)$-torus knots are $H$-restricted. The figure 
eight knot is $H$-restricted. 

\begin{prop} If $K_1$ and $K_2$ are $H$-restricted then $K_1\# K_2$ 
is $H$-restricted.
\end{prop}

\begin{conjecture} All knots are $H$-restricted. 
\end{conjecture} 

This is a homological counterpart of Conjecture 1 in [BN] about $Kh_{\Q}.$

\begin{prop} If $K$ is $H$-restricted then $\mathrm{rank}\cH(K)= 
   \mathrm{rank}\widetilde{\cH}(K)-1.$ 
\end{prop}

\section{Cohomology with $\Z_2$-coefficients} \label{Z2coefficients}

Let us now work over $\Z$ rather that $\Q,$ so that $\mc{A}= \Z[X]/(X^2).$ 
A computation in [Kh1, Section 6.2] implies that $\mc{C}(K),$ where $K$ is a
 $(2,2m+1)$-torus knot, is isomorphic to the direct sum (modulo
acyclic complexes) of the complex 
$0 \lra \mc{A}\{ 2m\}\lra 0 $ and $m$ complexes $C'_1$
\begin{equation}\label{integral}   
   0 \lra \mc{A} \stackrel{2X}{\lra} \mc{A}\{-2\} \lra 0
\end{equation} 
with various shifts.

Cohomology of $C'_1\otimes_{\Z}\Q$
is two-dimensional (over $\Q$), and is a matching pair of 
cohomology groups in bidegrees that differ by $(1,-4).$ 

\begin{figure} [htb] \drawing{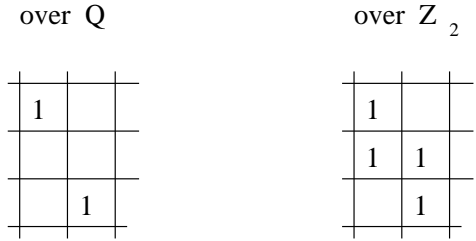} \caption{Dimensions of
cohomology of $C'_1$ over $\Q$ and $\Z_2$} \label{changedim} 
\end{figure}

Now change the base field to $\Z_2.$ In characteristic $2$ the 
differential in (\ref{integral}) is $0,$ and the dimension of cohomology 
groups of $C'_1\otimes_{\Z}\Z_2$ is $4$ (as a $\Z_2$-vector space), 
see figure~\ref{changedim}. 

According to the tables in Bar-Natan [BN], 
the same patterns relates rational and $\Z_2$-cohomology of any prime knot 
with at most seven crossings. Pair up the rational cohomology groups 
as in (ii), so that all but one pair look as on the 
left hand side of figure~\ref{changedim}, and change each on them 
to the quadruple of 1's on the right hand side. We get ranks of 
$\Z_2$-cohomology groups. 

It is likely that for any knot $K$ with at most 7 crossings
$\mc{C}(K)$ decomposes as a direct sum of 
\begin{itemize}
\item  
an acyclic complex, 
\item 
complex $0\lra \mc{A} \{ i\} \lra 0$ for some $i\in \Z,$ 
\item 
complexes $0 \lra \mc{A}\{j\} \stackrel{2kX}{\lra} \mc{A}\{j-2\} \lra
0$ for $j,k\in \Z.$ 
\end{itemize}

This would explain the observed relation between rational and 
$\Z_2$-cohomology of these knots.

\section{Cohomology of adequate knots} \label{adequate}

\begin{figure} [htb] \drawing{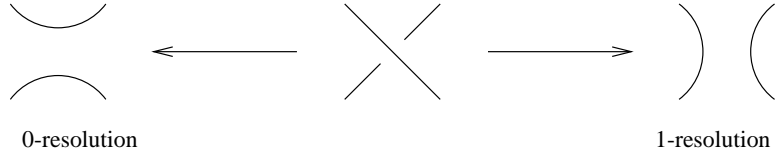} \caption{Two resolutions of a
crossing}  \label{resol} 
\end{figure}

For a link diagram $D$ denote by $s_+D$ (respectively $s_-D$) 
the diagram obtained by taking  0-resolution (respectively
1-resolution) of each crossing of $D,$ see figure~\ref{resol}.

We say that $D$ is \emph{adequate} if 
\begin{itemize}
\item for any crossing of $D$ the two segments of $s_+D$ that replace 
  this crossing belong to distinct components of $s_+D,$
\item  for any crossing of $D$ the two segments of $s_-D$ that replace 
  this crossing belong to distinct components of $s_-D.$
\end{itemize} 

A reduced alternating link diagram is adequate. A link admitting an 
adequate diagram is called \emph{adequate}. For further information about 
adequate links see Thistlethwaite [Th] and Lickorish [Li, Chapter 5]. 

\begin{prop} Adequate non-alternating knots are H-thick. 
\end{prop}

\emph{Proof: } Assume $D$ is an adequate non-alternating diagram of a
knot $K.$ We continue to use cohomology with integer coefficients. 
Recall from [Kh1, Chapter 7] that
\begin{equation*}
  \overline{\cH}^0(D) \not= 0 \not= \overline{\cH}^n(D)
\end{equation*}
where $n$ is the number of crossing of $D.$ More precisely, 
\begin{eqnarray} 
   \overline{\cH}^{0,|s_+D|}(D)\cong \Z, & \hspace{0.1in} & 
    \overline{\cH}^{0,i}(D)\cong 0, 
     \hspace{0.05in} \mathrm{if} \hspace{0.05in}
      i>|s_+D|, \\
   \overline{\cH}^{n,-|s_-D|-n}(D)\cong \Z, & \hspace{0.1in} & 
    \overline{\cH}^{0,i}(D)\cong 0, 
     \hspace{0.05in} \mathrm{if} \hspace{0.05in}
      i<-|s_-D|-n, 
 \end{eqnarray}
where $|s_+D|$ is the number of components of $s_+D,$ etc.

\begin{figure} [htb] \drawing{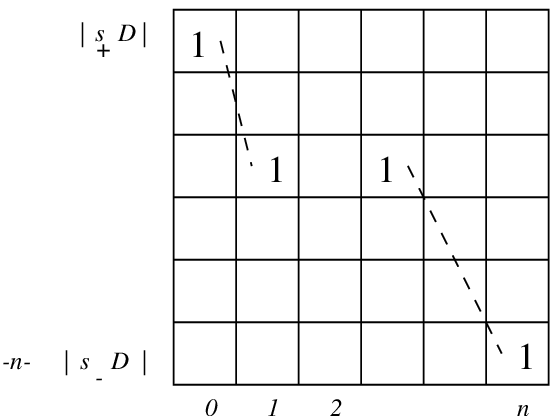} \caption{} \label{adeqdiag} 
\end{figure}

From discussion is Section~\ref{mod-structure} we know that rational cohomology
 groups come in 
pairs (complex (\ref{indec}) contributes $\Q\oplus \Q$ to cohomology, in 
two degrees that differ by  $(n,-2n-2)$). The companion of 
$\overline{\cH}^{0,|s_+D|}(D)\otimes \Q\cong \Q$ will lie one diagonal 
below it, while 
the companion of $\overline{\cH}^{n,-|s_-D|-n}\otimes \Q$ will lie one diagonal
above it. This is illustrated in figure~\ref{adeqdiag}, which
unintentionally shows the case $n = |s_+D| + |s_-D|.$ 

If $K$ is $H$-thin, these two pairs of cohomology groups must lie on
two adjancent diagonals. This implies $n+2 = |s_+D| + |s_-D|.$ 

\begin{lemma} If $D$ is adequate, non-alternating and prime then 
$n+2 > |s_+D| + |s_-D|.$ If $D$ is alternating then $n+2 = |s_+D| + |s_-D|.$ 
\end{lemma} 
 This lemma is proved in [Li, Chapter 5]. $\square$ 

Therefore, if $D$ is prime, $K$ is $H$-thick. The case of composite
$D$ follows from Corollary~\ref{con-sum}. 
$\square$

There are no adequate non-alternating knots with 9 or fewer crossings, 
3 adequate non-alternating knots with 10 crossings: 
$10_{152},$  $10_{153},$  $10_{154},$ and 15 adequate non-alternating 
11-crossing knots.

\section{Cohomology of positive and braid positive knots} 

\vspace{0.05in}

{\bf Positive knots}

\vspace{0.05in}

\begin{figure} [htb] \drawing{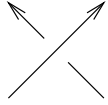} \caption{A positive crossing} 
 \label{positive} 
\end{figure}

We say that a knot is \emph{positive} if it has a diagram with only positive 
crossings (figure~\ref{positive}). 

\begin{prop} If $K$ is a positive knot then $\cH^{i,j}(K)=0$ if 
$i<0,$ 
   \[ \cH^{0,j}(K) = \left\{ \begin{array}{ll} 
                            \Z   &    \mbox{if $j = s-n-1 \pm 1$} \\
                             0   & \mbox{otherwise,}
                           \end{array}
                   \right.   \]
 and $\cH^{i,j}=0$ if $i>0$ and $j\ge s-n ,$ where $s$ is the number 
of Seifert circles and $n$ the number of crossings in a positive 
diagram of $K.$ 
\end{prop} 

\emph{Proof:} Left to the reader. $\square$   

Note that $\frac{n - s + 1}{2}$ is the genus of $K.$ 

\vspace{0.1in}

{\bf Braid positive knots}

\vspace{0.1in}

$8_{19}$ is a $(3,4)$-torus knot, $10_{124}$ is a $(3,5)$-torus 
knot. Both are H-thick. If $n,m$ are odd, the $(n,m)$-torus knot is
H-thick since its Jones polynomial is not alternating.  
We expect that $(n,m)$-torus knots, $2<n<m,$ are H-thick. 

Torus knots are examples of \emph{braid positive} knots, i.e. knots that are 
closures of positive braids. 

Braid positive prime knots with at most 10 crossings are $(2,n)$-torus knots, 
for $n\in \{ 3,5,7,9\},$ and the four H-thick knots $8_{19},$  $10_{124},$ 
 $10_{139},$  $10_{152}.$ 

\vspace{0.05in}

There are two braid positive prime 11-crossing knots: the
$(2,11)$-torus knot and $11^n_{77},$ the closure of the 
braid $\sigma_1^2 \sigma_2^2 \sigma_1 \sigma_3 \sigma_2^3 \sigma_3^2.$ 
The latter is H-thick [BN]. 

\vspace{0.05in}

There are 7 braid positive prime knots with 12
crossings. All of them are H-thick, since their Jones polynomials 
are not alternating. 

\vspace{0.05in} 

Not counting the $(2,13)$-torus knot, there are 12 braid positive 
prime 13-crossing knots. At least 10 are H-thick  (the Jones
polynomial is not alternating). We don't know if the remaining knots 
$13^n_{4587}$ and $13^n_{5016}$ are H-thick. 

\vspace{0.05in}

There are 17 braid positive prime knots with 14 crossings. All but 
3 have non-alternating Jones polynomial. 

\vspace{0.1in}

{\bf Problem:}  \emph{Are all braid positive prime knots other than 
  $(2,n)$-torus knots H-thick?} 
  
\vspace{0.1in}

{\bf Problem:} \emph{If $K$ is braid positive, is $\cH^{1,j}(K)=0$ for all
$j$?}

\section{Alexander polynomial and cohomology} \label{alexander}

We say that a prime 
knot is \emph{Ap-special} if its Alexander polynomial is not  
alternating or has a gap. 
A well-known theorem of Murasugi [Mu] can be restated as 

\begin{prop} Ap-special knots are not alternating. 
\end{prop} 

Few small knots are Ap-special, and all or nearly all small Ap-special
knots are H-thick: 

\begin{itemize} 
\item There are 9  
 Ap-special knots with at most 10 crossings. All of them are H-thick. 
\item  There are 19  Ap-special knots with 11
 crossings. All of them are H-thick.  
\item There are 104 Ap-special knots with 12 
 crossings. For all but 8 of them the Jones polynomial is not
alternating, so that at least 96 of these knots are H-thick. 
\item 
There are 115 knots with 13 crossings and a gap in the Alexander 
polynomial. All but 13 have non-alternating Jones polynomial, thus, 
at least 102 of these knots are H-thick. 
\end{itemize}

\vspace{0.05in}

{\bf Problem:} \emph{Is any Ap-special knot  H-thick?} 

\vspace{0.05in} 

Knots with non-alternating Jones polynomial are a minority 
among non-alternating knots with at most 14 crossings, as seen 
in the table below.

\[
\begin{array}{ccccccc}  
  \mathrm{crossings}          & \le 9 & 10 & 11  & 12  & 13   & 14    \\
  \mathrm{non-alternating}    &    11 & 42 & 185 & 888 & 5110 & 27110 \\
  \mathrm{Jones \hspace{0.1in} not \hspace{0.1in} alternating}     
                              &     0 &  7 & 26  & 169 & 1154 & 7075  \\
  \mathrm{H-thick}            &   2   & 10 & 41  & \ge 169 & \ge 1154
  &   \ge 7075
  \end{array}
\]

For instance, the fifth column says that there are 888 prime
non-alternating knots with 12 crossings (not distinguishing mirror
images); among them 169 have non-alternating Jones polynomial, and, 
therefore, at least 169 are H-thick. 
On the other hand, there is no doubt that for large $n$ most $n$-crossing
knots are H-thick.

\vspace{0.1in}

The following examples provide another experimental relation between 
the Alexander polynomial and knot cohomology. 

\begin{enumerate} 
\item The only knot with the trivial Alexander polynomial and at most 
10 crossings is the unknot. There are two 11-crossing,  two 12-crossing, 
 fifteen 13-crossing and thirty-six 14-crossing knots with the trivial 
 Alexander polynomial. All of them are H-thick (since their Jones
 polynomials are not alternating). 
\item The Alexander polynomial of the trefoil is $t^{-1}-1 +t.$ There 
are no other knots with at most 12 crossings and this Alexander
polynomial. There are eight 13-crossing knots and seventeen prime
14-crossing knots with this Alexander polynomial. All of them are 
H-thick (for the same reason). 
\item The figure eight knot is the only one with less than  13
crossings and  Alexander polynomial $-t^{-1}+3-t.$ There are two 
13-crossing knots and fifteen 14-crossing knots with this Alexander 
polynomial. All are H-thick. 
\item $\Delta(5_2)= 2t^{-1} - 3 + 2t.$ There are no other knots with 
this Alexander polynomial and less than $12$ crossings. Four 
12-crossing, three 13-crossing, and nine 14-crossing knots have 
Alexander polynomial $2t^{-1} - 3 + 2t.$ All of these knots are 
H-thick. 
\item Consider knots with at most 14 crossings and Alexander  
polynomial $-2t^{-1}+5-2t.$ Four of them: $6_1,$  $9_{46},$ 
 $11^n_{139},$ and $13^n_{3523}$  are H-thin (these knots are examples
 of $(n,-3,3)$-pretzel knots; any $(n,-3,3)$-pretzel knot 
is slice, H-thin, and its cohomology has rank $10$). 
The remaining two  11-crossing knots, 
four 12-crossing knots, eleven 13-crossing knots, and fifty
14-crossing knots with this polynomial are H-thick.   
\item $\Delta(5_1)=\Delta(10_{132})= t^{-2} - t^{-1} + 1 - t + t^2.$ 
$5_1$ is the $(2,5)$-torus knot and is H-thin. $10_{132}$ is H-thick.
There are no 11 and 12-crossing knots with this Alexander polynomial. 
Two 13-crossing knots and twelve 14-crossing knots have this Alexander 
polynomial. All are H-thick. 
\item $10_{153}$ is the only knot with at most 11 crossing and 
Alexander polynomial $t^{-3}- t^{-2} - t^{-1} + 1 - t - t^2 + t^3.$ 
Four 12-crossing, seven 13-crossing and and nineteen
14-crossing knots have this Alexander polynomial. All are H-thick. 
Unlike other examples, this Alexander polynomial is not
alternating. 
\end{enumerate}

These examples suggest that knots with small Alexander 
polynomial relative to the crossing number tend to be H-thick.

\section{Volume and cohomology} 

H-thick knots with few crossings tend to have small hyperbolic volume 
or to be non-hyperbolic: 
\begin{itemize}
\item  
$8_{19}$ is the only H-thick knot with 8 crossings and the only 
non-hyperbolic knot with 8 crossings (it is the $(3,4)$-torus knot). 
\item The H-thick knot $10_{124}$ is the $(3,5)$-torus knot and the 
only non-hyperbolic 10-crossing knot. 
\item 
$9_{42},$ the only H-thick 9-crossing knot, has the second 
smallest volume ($\approx 4.05686$) among all 48   
hyperbolic knots with 9 crossings (and
the smallest determinant ($=7$) among all 9-crossing knots).  
$9_{42}$ has the same volume as $10_{132}$, 
another H-thick knot. The latter has the smallest volume 
among all hyperbolic knots with 10 crossings. 
Among known pairs of knots with the same volume, $(9_{42},10_{132})$
is the pair with the 
second smallest volume. The pair with the smallest volume 
($\approx 2.8281$) consists of $5_2$ and 
the famous $(-2,3,7)$-pretzel knot. 
Knot $5_2$ is H-thin, while the $(-2,3,7)$-pretzel knot 
is H-thick, since its Jones polynomial is not alternating.  
\item Three out of the four hyperbolic 10-crossing knots with the smallest 
 volumes are H-thick, even though among 164 hyperbolic 10-crossing knots
only 9 are H-thick.   
\end{itemize}

\vspace{0.1in}

Determinant $det(K)$ of a knot $K$ is the determinant of the matrix 
$M+M^{T}$ where $M$ is a Seifert matrix of $K.$ Determinant is a 
common specialization of the Alexander and Jones polynomials:  
\begin{equation*}
  det(K) = \Delta_{-1}(K) = J_{\sqrt{-1}}(K)
\end{equation*} 
$|det(K)|$ is also the number of elements in the first homology 
group of the double cover of $\mathbb{S}^3$ branched over $K.$ 

\vspace{0.1in}

Nathan Dunfield  documents a fascinating relation between determinants
and volumes of hyperbolic knots [D]. First, he plots $\log|det(K)|$ 
versus the volume of $K$ for all alternating knots $K$  
with a fixed number of crossings. Amazingly, the points cluster around 
a straight line. Next, he combines the pictures into one by plotting 
$\frac{\log|det(K)|}{\log(\mathrm{deg}J(K))}$ versus the volume of $K$ for all
alternating knots with at most 13 crossings and samples of 14-16
crossing alternating knots. Again, all points stay close to a straight line.  

Dunfield comments: "$\log(J(-1))$ is one of the first terms in
Kashaev's conjecture about the relationship between the colored Jones 
polynomial and hyperbolic volume. However, the above doesn't appear 
to simply be saying that you have fast convergence in Kashaev's
conjecture as the slope of the line is not what you would expect."

\vspace{0.1in}

\begin{figure}    
\readdata{\junk}{vj10.dat}
\begin{pspicture}(1,-2.6)(15,5)
\psset{xunit=0.7cm, yunit=0.7cm}
\rput(8.5,-2.4){Volume($K$)}
\rput{90}(1.5,1.8){$\log|det(K)|$} 
\psaxes[axesstyle=frame,tickstyle=bottom,Ox=3,Oy=-1](3,-1)(3,-1)(14,4.5) 
\dataplot[plotstyle=dots]{\junk}
\end{pspicture}
\caption{Volume versus $\log|det(K)|$ for 10-crossing 
non-alternating knots}\label{det10}
\end{figure}
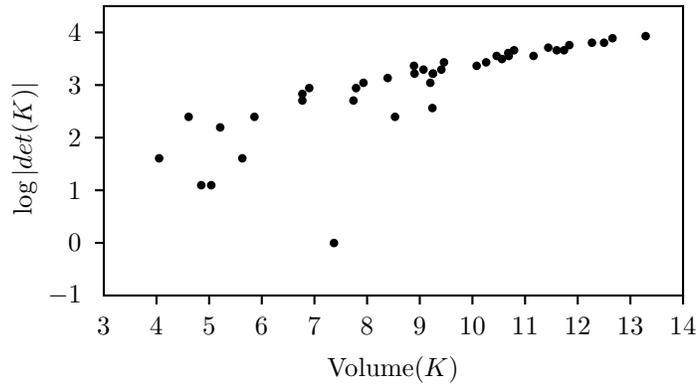

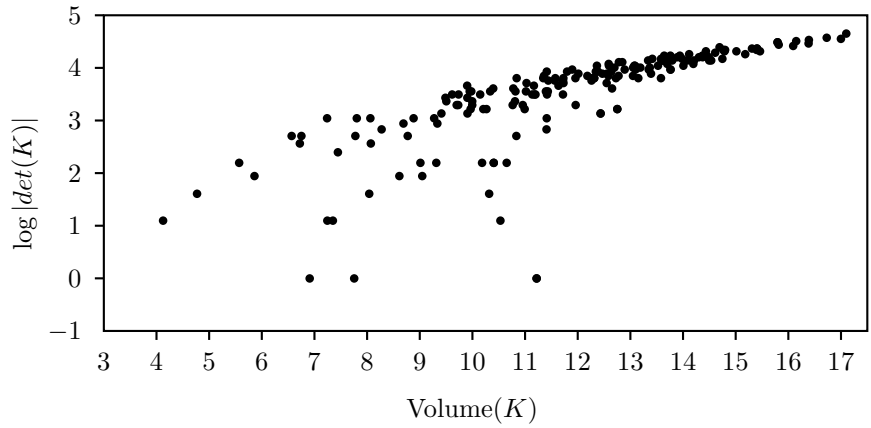
\begin{figure}
\readdata{\junk}{vj11.dat}
\begin{pspicture}(1,-2.6)(18,6)
\psset{xunit=0.7cm, yunit=0.7cm}
\rput(10,-2.5){Volume($K$)}
\rput{90}(1.5,1.8){$\log|det(K)|$} 
\psaxes[axesstyle=frame,tickstyle=bottom,Ox=3,Oy=-1](3,-1)(3,-1)(17.5,5) 
\dataplot[plotstyle=dots]{\junk}
\end{pspicture}
\caption{Volume versus $\log|det(K)|$ for 11-crossing 
non-alternating knots}\label{det11}
\end{figure}

When non-alternating knots are included, the plots become less 
impressive. The majority of points still lie close to the coveted 
straight line, but there are defections.  For instance, 
there are hyperbolic knots 
with $det(K)= \pm 1,$ and points assigned to them will lie on the
x-axis, far away from where we would like them to. This is 
illustrated in figures~\ref{det10}, \ref{det11}, 
where we plot $(vol(K),\log|det(K)|)$ for all hyperbolic
non-alternating knots with 10 and 11 crossings (for 12, 13 crossings 
consult [D]). 

\vspace{0.06in}

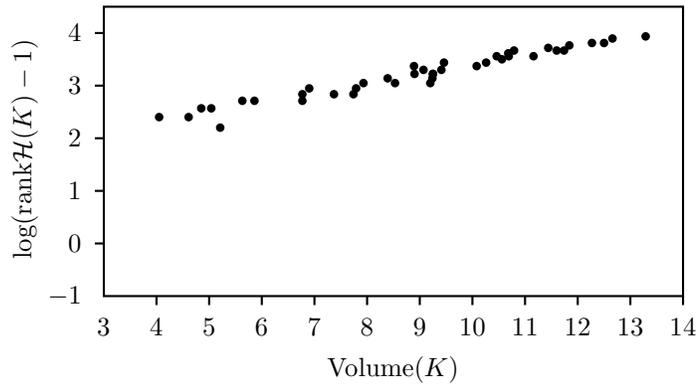
\begin{figure}   
\readdata{\junk}{vh10.dat}
\begin{pspicture}(1,-2.6)(15,5)
\psset{xunit=0.7cm, yunit=0.7cm}
\rput(8.5,-2.4){Volume($K$)}
\rput{90}(1.5,1.8){$\log(\mathrm{rank}\cH(K)-1)$} 
\psaxes[axesstyle=frame,tickstyle=bottom,Ox=3,Oy=-1](3,-1)(3,-1)(14,4.5) 
\dataplot[plotstyle=dots]{\junk}
\end{pspicture}
\caption{Volume versus $\log(\mathrm{rank}\cH(K)-1)$ for 10-crossing 
non-alternating knots}\label{homol10}
\end{figure}

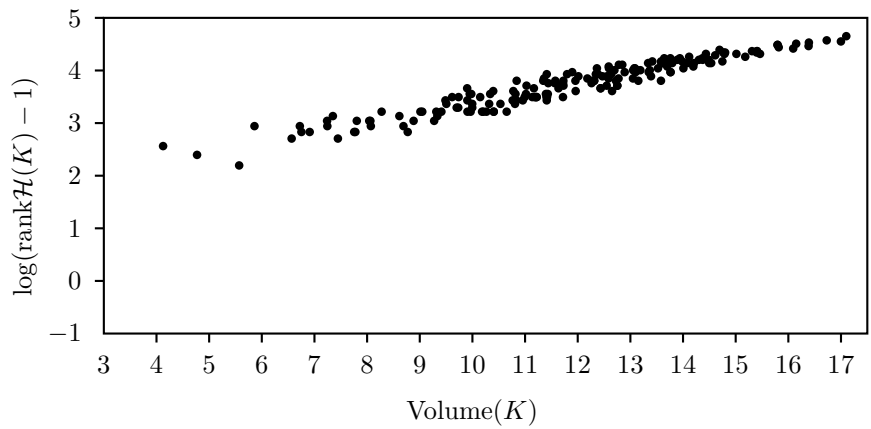
\begin{figure} 
\readdata{\junk}{vh11.dat}
\begin{pspicture}(1,-2.6)(18,6)
\psset{xunit=0.7cm, yunit=0.7cm}
\rput(10,-2.5){Volume($K$)}
\rput{90}(1.5,1.8){$\log(\mathrm{rank}\cH(K)-1)$} 
\psaxes[axesstyle=frame,tickstyle=bottom,Ox=3,Oy=-1](3,-1)(3,-1)(17.5,5) 
\dataplot[plotstyle=dots]{\junk}
\end{pspicture}
\caption{Volume versus $\log(\mathrm{rank}\cH(K)-1)$ for 11-crossing 
non-alternating knots}\label{homol11}
\end{figure}

To save the day, we change from $det(K)$ to the rank of the reduced
cohomology group of $K.$ The inequality 
\begin{equation*} 
   \mathrm{rank}\widetilde{\cH}(K) \ge |det(K)| 
\end{equation*} 
is valid for all knots, and turns into equality for H-thin knots. 
If the knot is $H$-restricted (and we expect that all knots are), 
$\mathrm{rank}\widetilde{\cH}(K)=\mathrm{rank} \cH(K)-1 .$ 
In figures~\ref{homol10}, \ref{homol11} 
we plot $(vol(K),\log(\mathrm{rank}\cH(K)-1))$ for all 
hyperbolic non-alternating knots with 10 and 11 crossings (there 
are 41, respectively 185, such knots).   

Clearly, for non-alternating knots with 10 and 11 crossings the 
correlation between the volume and the rank of cohomology is even better 
than the one between the volume and the determinant. 
Somehow $vol(K)$ and $\mathrm{rank}\cH(K)$ are successful in spying 
on each other. We have no explanation for this behaviour.

\clearpage

\clearpage

\clearpage

\section{Acknowledgments}

This paper owes its existence to Dror Bar-Natan and his work [BN]. I
greatly benefitted from 
discussions with Dror Bar-Natan, Michael Hutchings, Vaughan Jones, 
Greg Kuperberg, Paul Seidel, Alexander Stoimenow, and many others. 
Greg Kuperberg skillfully guided me through the jungles of PStricks, 
Alexander Stoimenow provided the list of braid positive knots 
with 14 or fewer crossings. \emph{Knotscape} was put to heavy use, and I am 
very much indebted to its creators Jim Hoste,  Morwen  Thistlethwaite, 
and Jeffrey Weeks.

\section{References}

[BN] D. Bar-Natan, On Khovanov's categorification of the Jones 
polynomial, arXiv:math.QA/0201043.    

[D] N. Dunfield, Jones polynomial and hyperbolic volume, \newline 
 www.math.harvard.edu/$\sim$nathand. 

[G] S. Garoufalidis, A conjecture on Khovanov's invariants, 
 University of Warwick preprint, October 2001.  

[Kh1] M. Khovanov, A categorification of the Jones polynomial, 
 \emph{Duke Math J.}, 101 (3), 359-426, 1999, 
 arXiv:math.QA/9908171.  

[Kh2] M. Khovanov, A functor-valued invariant of tangles, \newline
 arXiv:math.QA/0103190. 

[Lee] E. S. Lee, The support of the Khovanov's invariants for alternating 
  knots, arXiv:math.GT/0201105.   

[Li] W.B.R. Lickorish, An introduction to knot theory, 
 Graduate texts in mathematics 175, 1997, Springer-Verlag. 

[Mu] K. Murasugi, On the Alexander polynomial of the alternating 
 knot, Osaka Math. J., 10 1958, 181-189; errata, 11 (1959), 95. 

[St] A. Stoimenow, Polynomials of knots with up to 10 crossings,  
 \newline 
 http://guests.mpim-bonn.mpg.de/alex. 

[Th] M.B. Thistlethwaite, On the Kauffman polynomial of an adequate 
 link, Invent. Math., {\bf 93}, (1988), 285-296.

\end{document}